\documentclass[11pt,a4paper]{article}

\usepackage[margin=27mm]{geometry}
\usepackage{amsmath,amssymb,amsthm,mathtools,mathrsfs}
\usepackage{microtype}
\usepackage{needspace}
\usepackage[hidelinks,unicode]{hyperref}
\usepackage{aliascnt}
\usepackage[nameinlink,noabbrev]{cleveref}
\usepackage{indentfirst}

\hypersetup{
  pdftitle={Every Extended Metric Space Is Isometrically Embeddable into the Gromov-Hausdorff Class},
  pdfauthor={Ryo Fukano}
}

\allowdisplaybreaks[2]

\theoremstyle{plain}
\newtheorem{theorem}{Theorem}[section]

\newaliascnt{lemma}{theorem}
\newtheorem{lemma}[lemma]{Lemma}
\aliascntresetthe{lemma}

\newaliascnt{proposition}{theorem}
\newtheorem{proposition}[proposition]{Proposition}
\aliascntresetthe{proposition}

\newaliascnt{corollary}{theorem}
\newtheorem{corollary}[corollary]{Corollary}
\aliascntresetthe{corollary}

\theoremstyle{definition}
\newaliascnt{definition}{theorem}
\newtheorem{definition}[definition]{Definition}
\aliascntresetthe{definition}

\newaliascnt{example}{theorem}

\aliascntresetthe{example}

\theoremstyle{remark}
\newaliascnt{remark}{theorem}

\aliascntresetthe{remark}

\crefname{theorem}{theorem}{theorems}
\crefname{lemma}{lemma}{lemmas}
\crefname{proposition}{proposition}{propositions}
\crefname{corollary}{corollary}{corollaries}
\crefname{definition}{definition}{definitions}
\crefname{example}{example}{examples}
\crefname{remark}{remark}{remarks}
\crefname{section}{section}{sections}
\crefname{subsection}{subsection}{subsections}
\crefname{equation}{equation}{equations}
\Crefname{theorem}{Theorem}{Theorems}
\Crefname{lemma}{Lemma}{Lemmas}
\Crefname{proposition}{Proposition}{Propositions}
\Crefname{corollary}{Corollary}{Corollaries}
\Crefname{definition}{Definition}{Definitions}
\Crefname{example}{Example}{Examples}
\Crefname{remark}{Remark}{Remarks}
\Crefname{section}{Section}{Sections}
\Crefname{subsection}{Subsection}{Subsections}
\Crefname{equation}{Equation}{Equations}

\newcommand{\M}{\mathcal{M}}
\newcommand{\GH}{\mathfrak{GH}}
\newcommand{\R}{\mathbb R}
\newcommand{\N}{\mathbb N}
\newcommand{\dGH}{d_{\mathrm{GH}}}
\newcommand{\dH}{d_{\mathrm H}}
\newcommand{\dis}{\operatorname{dis}}
\newcommand{\diam}{\operatorname{diam}}
\newcommand{\dens}{\operatorname{dens}}
\newcommand{\Corr}{\operatorname{Corr}}
\newcommand{\Aut}{\operatorname{Aut}}

\newcommand{\id}{\operatorname{id}}

\title{Every Extended Metric Space Is Isometrically Embeddable into the Gromov--Hausdorff Class}
\author{Ryo Fukano}
\date{October 27, 2026}

\begin{document}
\maketitle
\begin{NoHyper}
\begingroup
\renewcommand{\thefootnote}{}
\endgroup
\end{NoHyper}

\begin{abstract}
Let $\GH$ denote the Gromov--Hausdorff class of isometry classes of nonempty metric spaces, equipped with the Gromov--Hausdorff extended pseudometric.
We prove that every nonempty extended metric space is isometrically embeddable into $\GH$.
More precisely, for every nonempty extended metric space $X$, all image points may be represented by metrics on a common set of cardinality $\max\{\dens X,\aleph_0\}$.
The main ingredient is a realization theorem for function spaces: for every nonempty set $I$, the space of all real-valued functions on $I$, endowed with the uniform extended metric, can be realized as a family of metrics on a fixed underlying set. Under this realization, the Gromov--Hausdorff distance agrees exactly with the uniform extended distance, including the case of infinite distance. The resulting metric spaces are complete, uniformly discrete, and unbounded, and have cardinality $\max\{|I|,\aleph_0\}$. If $I$ is countable, they may moreover be chosen proper and of bounded geometry.
As further consequences, we obtain sharp cardinality results for clouds and universality results for cardinality-restricted parts of $\GH$ and for spheres in the associated zero-distance quotients.

\end{abstract}

\section{Introduction}

Let $\M$ denote the space of isometry classes of nonempty compact metric spaces equipped with the Gromov--Hausdorff distance, and let $\GH$ denote the Gromov--Hausdorff class of isometry classes of all nonempty metric spaces. The Gromov--Hausdorff distance is a finite-valued metric on $\M$, whereas on $\GH$ it is an extended pseudometric: distances may be infinite, and nonisometric spaces may have distance zero.

The question of whether every nonempty compact metric space admits an isometric embedding into $\M$ was posed by S.~Iliadis in \cite{IIT}, where the authors established an affirmative answer for finite metric spaces. Byakuno restated the question in \cite[Problem~1.1, Theorem~1.3]{Byakuno} and proved an isometric embedding theorem for countable products with summable diameter bounds. In \cite{FukanoCompact}, the author proved the assertion for every nonempty compact metric space. For every $D>0$ and every nonempty compact metric space $K$ with $\diam K\le D$, the proof assigns to each $[0,D]$-valued $1$-Lipschitz function on $K$ a metric on a fixed Cantor space, with Gromov--Hausdorff distances equal to the uniform distances between functions.

For the target $\GH$, Ivanov--Tuzhilin \cite[Theorem~3]{ITclass} proved that every bounded metric space admits an isometric embedding into $\GH$, leaving the unbounded case open. Abdullaev \cite[Theorem~4.4]{Abdullaev} subsequently constructed metric spaces whose Gromov--Hausdorff distances agree with the Hausdorff distances between subsets of $[0,\infty)$ containing $0$. From this construction he deduced that $\ell^\infty(\N)$ and hence every separable metric space, including unbounded ones, admit isometric embeddings into $\GH$ \cite[Theorem~4.6, Corollary~4.7]{Abdullaev}. This construction also allows infinite distances, but these results do not establish isometric embeddability for arbitrary nonseparable metric spaces or arbitrary extended metric spaces.

Our main results are the following. The notation is defined in \cref{sec:preliminaries}.

\Needspace{6\baselineskip}
\begin{theorem}\label{thm:main}
For every nonempty extended metric space $(X,d_X)$, set $\kappa:=\max\{\dens X,\aleph_0\}$. Then there exists an isometric embedding
\[
 \Phi:(X,d_X)\longrightarrow(\GH_{\le\kappa},\dGH).
\]
Moreover, one can choose a set $\Omega$ and a family $(\rho_x)_{x\in X}$ of metrics on $\Omega$ such that $\Phi(x)=[(\Omega,\rho_x)]$ for every $x\in X$. Each $(\Omega,\rho_x)$ is complete, uniformly discrete, and unbounded, and
\begin{equation}\label{eq:main}
 \begin{gathered}
 |\Omega|=\kappa,\\
 \dGH\bigl((\Omega,\rho_x),(\Omega,\rho_y)\bigr)=d_X(x,y)
 \quad(x,y\in X).
 \end{gathered}
\end{equation}
If $d_X(x,y)<\infty$, the identity relation $I_\Omega:=\{(\omega,\omega)\mid\omega\in\Omega\}$ is an optimal correspondence between $(\Omega,\rho_x)$ and $(\Omega,\rho_y)$.
\end{theorem}

The main theorem follows from the following realization theorem for a function space.

\Needspace{6\baselineskip}
\begin{theorem}\label{thm:coding-intro}\label{thm:coding}
Let $I$ be a nonempty set and let $\Delta_I$ be the uniform extended metric on $\R^I$ defined in \cref{def:uniform-function-space}. There exist a set $\Omega_I$, independent of $f$, and a finite-valued metric $\rho_f$ on $\Omega_I$ for every $f\in\R^I$, such that, for all $f,g\in\R^I$,
\begin{equation}\label{eq:coding}
 \begin{aligned}
 \dGH\bigl((\Omega_I,\rho_f),(\Omega_I,\rho_g)\bigr)
 &=\Delta_I(f,g)\\
 &=\frac12\sup_{\omega,\omega'\in\Omega_I}
       |\rho_f(\omega,\omega')-\rho_g(\omega,\omega')|.
 \end{aligned}
\end{equation}
Each $(\Omega_I,\rho_f)$ is complete, uniformly discrete, and unbounded, and we may take $|\Omega_I|=\kappa_I:=\max\{|I|,\aleph_0\}$. Consequently, $f\mapsto[(\Omega_I,\rho_f)]$ defines an isometric embedding from $(\R^I,\Delta_I)$ into $(\GH_{\le\kappa_I},\dGH)$. In particular, if $\Delta_I(f,g)<\infty$, the identity relation $I_{\Omega_I}$ is an optimal correspondence.
\end{theorem}

This paper is organized as follows. \Cref{sec:preliminaries} collects the preliminaries on extended metric spaces, the Gromov--Hausdorff distance and correspondences, and clouds. In \cref{sec:construction}, we construct metrics on a common underlying set, prove \cref{thm:coding}, and deduce \cref{thm:main} using truncated distance functions. We also give proper representatives of bounded geometry when $I$ is countable. Finally, \cref{sec:clouds} establishes universality results for cardinality-restricted clouds and examines spheres in their zero-distance quotients.

\section{Preliminaries}
\label{sec:preliminaries}

Throughout, we assume the axiom of choice and work with classes in the sense of von Neumann--Bernays--G\"odel class theory. Unless otherwise stated, metric spaces are nonempty and the distance between any two points is finite.

\begin{definition}\label{def:extended-metric}
An \emph{extended metric} on a set $X$ is a map $d:X\times X\to[0,\infty]$ satisfying symmetry, the triangle inequality in the extended-real sense, and $d(x,y)=0\iff x=y$. If the last condition is weakened to $d(x,x)=0$, we call $d$ an \emph{extended pseudometric}. We use the same axioms for extended metrics and extended pseudometrics on classes.

The topology of an extended metric space $(X,d)$ is defined by taking the open balls $B_d(x,r):=\{y\in X\mid d(x,y)<r\}$ of finite radius $0<r<\infty$ as a base. Its \emph{density character} is
\[
 \dens X:=\min\{|D|\mid D\subseteq X\text{ is dense}\}.
\]
The equivalence classes of $x\sim y\iff d(x,y)<\infty$ are called \emph{finite-distance components}. We write $\mathcal A(X)$ for their index set and $(X_\alpha)_{\alpha\in\mathcal A(X)}$ for the components.
\end{definition}

\begin{definition}\label{def:isometric-universality}
Let $(X,d_X)$ be an extended metric space and $(Y,d_Y)$ a class equipped with an extended pseudometric. An injective map $F:X\to Y$ is an \emph{isometric embedding} if $d_Y(F(x),F(x'))=d_X(x,x')$ for all $x,x'\in X$. We say that $Y$ is \emph{isometrically universal} for a family $\mathscr A$ of metric or extended metric spaces if every member of $\mathscr A$ admits an isometric embedding into $Y$.
\end{definition}

\begin{definition}\label{def:uniform-function-space}
For a nonempty set $I$, let $\R^I$ denote the set of all real-valued functions on $I$, regarded as families $(f(i))_{i\in I}$ with \emph{index set} $I$. Define the \emph{uniform extended metric} on $\R^I$ by
\begin{equation}\label{eq:Delta}
 \Delta_I(f,g):=\sup_{i\in I}|f(i)-g(i)|\in[0,\infty]
 \qquad(f,g\in\R^I).
\end{equation}
Write
\[
 \ell^\infty(I):=\{h\in\R^I\mid \sup_{i\in I}|h(i)|<\infty\},\qquad
 \|h\|_\infty:=\sup_{i\in I}|h(i)|
\]
for the subspace of bounded functions and its uniform norm. Two functions $f,g\in\R^I$ represent the same coset in $\R^I/\ell^\infty(I)$ if and only if $\Delta_I(f,g)<\infty$.
\end{definition}

The function $\Delta_I$ is an extended metric on $\R^I$.

For a point $x$ in a metric space $(X,d)$, a real number $r\ge0$, and a nonempty subset $A\subseteq X$, set
\[
 \begin{aligned}
 \overline B_d(x,r)&:=\{y\in X\mid d(x,y)\le r\},&
 S_d(x,r)&:=\{y\in X\mid d(x,y)=r\},\\
 d(x,A)&:=\inf_{a\in A}d(x,a),&
 \diam_d(A)&:=\sup_{a,a'\in A}d(a,a').
 \end{aligned}
\]
We omit the subscript when the metric is clear.

\begin{definition}\label{def:representative-geometry}\label{def:bounded-geometry}
A subset $A$ of a metric space is \emph{$\delta$-separated}, for $\delta>0$, if the distance between every two distinct points of $A$ is at least $\delta$. A metric space is \emph{uniformly discrete} if it is $\delta$-separated for some $\delta>0$, and it is \emph{proper} if every closed bounded subset is compact. A uniformly discrete metric space $(Z,d_Z)$ has \emph{bounded geometry} if $\sup_{z\in Z}|\overline B_{d_Z}(z,R)|<\infty$ for every $R>0$.
\end{definition}

\begin{lemma}\label{lem:component-density}
Every finite-distance component $X_\alpha$ of a nonempty extended metric space $X$ is both open and closed, and the restriction of the distance to $X_\alpha$ is finite-valued. A subset $D\subseteq X$ is dense in $X$ if and only if $D\cap X_\alpha$ is dense in $X_\alpha$ for every $\alpha\in\mathcal A(X)$, and
\begin{equation}\label{eq:density-components}
 \dens X=\sum_{\alpha\in\mathcal A(X)}\dens X_\alpha.
\end{equation}
The sum on the right is a cardinal sum.
\end{lemma}

\begin{proof}
Every open ball of finite radius is contained in a single finite-distance component. Hence each component and its complement are open, which also proves the equivalence concerning density. If $D$ is dense, then
\[
 |D|=\sum_\alpha|D\cap X_\alpha|
 \ge \sum_\alpha\dens X_\alpha.
\]
Conversely, choose a dense subset $D_\alpha\subseteq X_\alpha$ of cardinality $\dens X_\alpha$ for each $\alpha$. The union $\bigcup_\alpha D_\alpha$ is dense in $X$ and has cardinality equal to the sum on the right.
\end{proof}

For the Gromov--Hausdorff distance and correspondences, we follow \cite[Section~1]{ITclass}.

\begin{definition}\label{def:GH-distance}
For nonempty subsets $A,B$ of a metric space $(Z,d_Z)$, define their Hausdorff distance by
\begin{equation}\label{eq:Hausdorff}
 \dH^Z(A,B):=\max\left\{
 \sup_{a\in A}d_Z(a,B),\quad \sup_{b\in B}d_Z(b,A)
 \right\}\in[0,\infty].
\end{equation}
For metric spaces $X,Y$, a triple $(X',Y',Z)$, where $Z$ is a metric space and $X',Y'$ are subsets of $Z$ isometric to $X,Y$, respectively, is called a \emph{realization} of $(X,Y)$. Define the Gromov--Hausdorff distance by
\begin{equation}\label{eq:GH-definition}
 \dGH(X,Y):=\inf\left\{
 r\in[0,\infty)\;\middle|\;
 \begin{array}{l}
 \text{there exists a realization }(X',Y',Z)\text{ of }(X,Y)\\
 \text{such that }\dH^Z(X',Y')\le r
 \end{array}
 \right\},
\end{equation}
where $\inf\varnothing=\infty$.
\end{definition}

\begin{definition}\label{def:GH-class}
The class of isometry classes of nonempty metric spaces is called the \emph{Gromov--Hausdorff class} and is denoted by $\GH$. We denote the element determined by a metric space $X$ by $[X]$ and set $\dGH([X],[Y]):=\dGH(X,Y)$. Let $\M$ denote the space of isometry classes of compact metric spaces. For an infinite cardinal $\kappa$, set
\begin{equation}\label{eq:cardinal-slice}
 \GH_{\le\kappa}:=\{[X]\in\GH\mid |X|\le\kappa\},
\end{equation}
with the restriction of $\dGH$.
\end{definition}

\begin{proposition}\label{prop:GH-set-representation}
The isometry class of each nonempty metric space can be represented by a set, and the collection of these representations forms a class. Moreover, for every infinite cardinal $\kappa$, $\GH_{\le\kappa}$ is a set and $|\GH_{\le\kappa}|\le 2^\kappa$.
\end{proposition}

\begin{proof}
Identify each nonzero cardinal $\lambda$ with its initial ordinal, and let $\operatorname{Met}(\lambda)$ denote the set of all finite-valued metrics on $\lambda$. The permutation group $\operatorname{Sym}(\lambda)$ acts on $\operatorname{Met}(\lambda)$ by
\[
 (\sigma\cdot d)(\alpha,\beta)
 :=d\bigl(\sigma^{-1}(\alpha),\sigma^{-1}(\beta)\bigr)
 \qquad(\alpha,\beta\in\lambda).
\]
Writing $\mathcal O_\lambda(d)$ for the orbit of $d\in\operatorname{Met}(\lambda)$, the isometry class of $(\lambda,d)$ is represented by the set $(\lambda,\mathcal O_\lambda(d))$. Indeed, two metric spaces are isometric if and only if their cardinalities agree and the corresponding metrics belong to the same orbit. These pairs form a class, and in the following we represent $[X]$ by such a pair.

Fix an infinite cardinal $\kappa$. For every cardinal $\lambda$ with $0<\lambda\le\kappa$,
\[
 |\operatorname{Met}(\lambda)|
 \le |[0,\infty)^{\lambda\times\lambda}|
 \le (2^{\aleph_0})^\kappa=2^\kappa.
\]
The set of cardinals $\lambda$ with $0<\lambda\le\kappa$ has cardinality at most $\kappa$. Hence
\[
 \GH_{\le\kappa}
 =\bigl\{(\lambda,\mathcal O_\lambda(d))\bigm|
 0<\lambda\le\kappa,\ d\in\operatorname{Met}(\lambda)\bigr\}
\]
is a set, and
\[
 |\GH_{\le\kappa}|
 \le\sum_{\substack{\lambda\text{ a cardinal}\\0<\lambda\le\kappa}}
 |\operatorname{Met}(\lambda)|
 \le\kappa\cdot2^\kappa=2^\kappa.
\]
\end{proof}

\begin{definition}\label{def:correspondence}
A \emph{correspondence} between nonempty sets $X,Y$ is a relation $R\subseteq X\times Y$ satisfying
\[
 \{x\in X\mid\exists y\in Y,\ (x,y)\in R\}=X,\qquad
 \{y\in Y\mid\exists x\in X,\ (x,y)\in R\}=Y.
\]
Write $\Corr(X,Y)$ for the set of correspondences. For subsets $A\subseteq X$, $B\subseteq Y$, and a point $x\in X$, write
\[
 \begin{aligned}
 R[A]&:=\{y\in Y\mid\exists x\in A,\ (x,y)\in R\},\\
 R^{-1}[B]&:=\{x\in X\mid\exists y\in B,\ (x,y)\in R\},\qquad
 R(x):=R[\{x\}].
 \end{aligned}
\]
The inverse correspondence is $R^{-1}:=\{(y,x)\mid(x,y)\in R\}$. If $X,Y$ have the same underlying set $S$, write $I_S:=\{(s,s)\mid s\in S\}$ for the identity relation.

For metric spaces $(X,d_X)$ and $(Y,d_Y)$, the \emph{distortion} of a correspondence $R$ is
\begin{equation}\label{eq:distortion}
 \dis R:=\sup\bigl\{|d_X(x,x')-d_Y(y,y')|
 \mid (x,y),(x',y')\in R\bigr\}\in[0,\infty].
\end{equation}
When $\dGH(X,Y)<\infty$, a correspondence $R$ satisfying $\dis R=2\dGH(X,Y)$ is called an \emph{optimal correspondence}.
\end{definition}

\begin{proposition}[{\cite[Assertion~1]{ITclass}}]\label{prop:correspondence}
For any metric spaces $A,B$,
\begin{equation}\label{eq:correspondence-formula}
 2\dGH(A,B)=\inf_{R\in\Corr(A,B)}\dis R.
\end{equation}
\end{proposition}

\begin{proposition}[{\cite[Section~1]{ITclass}}]\label{prop:GH-pseudometric}
The expression $\dGH([X],[Y]):=\dGH(X,Y)$ is independent of the choice of representatives and defines an extended pseudometric on $\GH$. That is, for all $P,Q,T\in\GH$,
\[
 \begin{gathered}
 \dGH(P,P)=0,\qquad \dGH(P,Q)=\dGH(Q,P),\\
 \dGH(P,T)\le\dGH(P,Q)+\dGH(Q,T).
 \end{gathered}
\]
However, $\dGH(P,Q)=0$ does not in general imply $P=Q$.
\end{proposition}

\begin{definition}\label{def:cloud}
Following \cite[Introduction]{BT}, the \emph{cloud} of a metric space $M$ is
\begin{equation}\label{eq:cloud}
 \mathcal C(M):=\{[A]\in\GH\mid\dGH(A,M)<\infty\}.
\end{equation}
For $P=[M]\in\GH$, we also write $\mathcal C(P):=\mathcal C(M)$.

For an infinite cardinal $\kappa$ and a metric space $M$ with $|M|\le\kappa$, write
\begin{equation}\label{eq:cloud-slice}
 \mathcal C_{\le\kappa}(M):=\mathcal C(M)\cap\GH_{\le\kappa}.
\end{equation}
Define the \emph{zero-distance quotient} by
\begin{equation}\label{eq:metric-cloud-quotient}
 \mathcal Q_\kappa(M)
 :=\mathcal C_{\le\kappa}(M)/\!\sim_0,
 \qquad
 [A]\sim_0[B]\iff\dGH(A,B)=0,
\end{equation}
and denote the metric induced on the quotient by $\dGH$ as well.
\end{definition}

By \cref{prop:GH-set-representation}, $\mathcal C_{\le\kappa}(M)$ and its quotient are sets. Moreover, the triangle inequality in \cref{prop:GH-pseudometric} shows that replacing a representative by another at zero distance does not change distances and that
\[
 \dGH(A,B)\le\dGH(A,M)+\dGH(M,B)<\infty
 \qquad([A],[B]\in\mathcal C_{\le\kappa}(M)).
\]
Thus $\mathcal Q_\kappa(M)$ is a finite-valued metric space.

\section{Construction of isometric embeddings}
\label{sec:construction}

\subsection{Construction of the metrics}

In this subsection we prove \cref{thm:coding}. Fix a nonempty set $I$.

\begin{definition}\label{def:rigidity}
Let $H=(V,E)$ be a simple undirected graph, and write $\Aut(H)$ for its automorphism group. The graph $H$ is \emph{rigid} if $\Aut(H)=\{\id_V\}$. For distinct $u,v\in V$, define the \emph{adjacency indicator} of $H$ by
\begin{equation}\label{eq:adjacency}
 \varepsilon_H(u,v)
 :=\begin{cases}
  1,&\{u,v\}\in E,\\
  0,&\{u,v\}\notin E.
 \end{cases}
\end{equation}
\end{definition}

\begin{proposition}[Shramov; {\cite[Construction~1]{ITclass}}]\label{prop:rigid-graph-existence}
For every infinite cardinal $\kappa$, there exists a rigid simple undirected graph $H=(V,E)$ with $|V|=\kappa$.
\end{proposition}

\begin{corollary}\label{cor:rigid-labels}
Let $I$ be a nonempty set, let $o\notin I$, and set $\kappa:=\max\{|I|,\aleph_0\}$. Then one can choose a simple undirected graph $H=(V,E)$ such that
\begin{equation}\label{eq:rigidity}
 I\sqcup\{o\}\subseteq V,
 \qquad |V|=\kappa,
 \qquad \Aut(H)=\{\id_V\}.
\end{equation}
Moreover, for every nonidentity bijection $\sigma:V\to V$, there exist distinct $u,v\in V$ such that
\begin{equation}\label{eq:adjacency-change}
 |\varepsilon_H(u,v)-\varepsilon_H(\sigma u,\sigma v)|=1.
\end{equation}
\end{corollary}

\begin{proof}
Since $\kappa$ is an infinite cardinal and $|I\sqcup\{o\}|\le\kappa$, there exists a set $V$ containing $I\sqcup\{o\}$ with $|V|=\kappa$. By \cref{prop:rigid-graph-existence}, take a rigid graph $(W,F)$ on a vertex set $W$ of cardinality $\kappa$. Choose a bijection $\tau:W\to V$ and put $E:=\bigl\{\{\tau(u),\tau(v)\}:\{u,v\}\in F\bigr\}$. Since $\tau$ is a graph isomorphism from $(W,F)$ onto $H:=(V,E)$, we have $\Aut(H)=\{\id_V\}$. If a bijection $\sigma:V\to V$ does not satisfy \eqref{eq:adjacency-change}, then $\varepsilon_H(u,v)=\varepsilon_H(\sigma u,\sigma v)$ for every pair of distinct vertices $u,v$. Thus $\sigma$ is an automorphism of $H$, and consequently $\sigma=\id_V$.
\end{proof}

For the remainder of this subsection, fix one graph $H=(V,E)$ as in \cref{cor:rigid-labels}.

\begin{definition}\label{def:multiscale-metric}
For $n\ge1$, define a map $c_n:\R\to\R$ and a positive number $a_n$ by
\begin{equation}\label{eq:clipping-scale}
 c_n(t):=\max\{-n,\min\{t,n\}\},
 \qquad
 a_n:=100^n.
\end{equation}
For $f\in\R^I$ and distinct $u,v\in V$, set
\begin{equation}\label{eq:perturbation}
 b_{n,f}(u,v)
 :=\begin{cases}
  2\bigl(n+c_n(f(i))\bigr),
   &\text{if }\{u,v\}=\{o,i\}\text{ for some }i\in I,\\
  0,&\text{otherwise}.
 \end{cases}
\end{equation}

Define the $n$th block and the underlying set by
\begin{equation}\label{eq:blocks}
 B_n:=\{n\}\times V,
 \qquad
 \Omega_I:=\bigsqcup_{n\ge1}B_n=\N_{\ge1}\times V.
\end{equation}
For each $f\in\R^I$, define a symmetric function $\rho_f$ to be zero on the diagonal and, for distinct points, by
\begin{align}
 \rho_f\bigl((n,u),(n,v)\bigr)
 &=3a_n+a_n\varepsilon_H(u,v)+b_{n,f}(u,v),
 &&u\ne v,
 \label{eq:internal-distance}\\
 \rho_f\bigl((n,u),(m,v)\bigr)
 &=10\max\{a_n,a_m\},
 &&n\ne m.
 \label{eq:external-distance}
\end{align}
Write $M_f:=(\Omega_I,\rho_f)$.
\end{definition}

\begin{proposition}\label{prop:metric}
For every $f\in\R^I$, the function $\rho_f$ is a finite-valued metric on $\Omega_I$, and for distinct $p,q\in B_n$,
\begin{equation}\label{eq:internal-range}
 3a_n\le \rho_f(p,q)\le5a_n.
\end{equation}
Furthermore, $M_f$ has cardinality $\kappa=\max\{|I|,\aleph_0\}$ and is complete, uniformly discrete, and unbounded.
\end{proposition}

\begin{proof}
By the definition of $c_n$, for all $s,t\in\R$,
\begin{equation}\label{eq:clipping-properties}
 |c_n(s)-c_n(t)|\le|s-t|,
 \qquad
 c_n(t)=t\quad(n\ge|t|).
\end{equation}
Also,
\begin{equation}\label{eq:perturbation-bound}
 0\le b_{n,f}(u,v)\le4n\le a_n,
 \qquad a_{n+1}=100a_n.
\end{equation}
This proves \eqref{eq:internal-range}.

We verify the triangle inequality. Take distinct points $p,q,r\in\Omega_I$.

First, if $p,q,r\in B_n$, then \eqref{eq:internal-range} gives
\[
 \rho_f(p,r)\le5a_n\le6a_n
 \le\rho_f(p,q)+\rho_f(q,r).
\]
The same estimate holds under any permutation of $p,q,r$.

Next, if $p,q\in B_n$ and $r\in B_m$ with $n\ne m$, put $D:=10\max\{a_n,a_m\}$. Then
\[
 \rho_f(p,r)=\rho_f(q,r)=D,
 \qquad 0<\rho_f(p,q)\le5a_n\le D.
\]
The three side lengths are $\rho_f(p,q),D,D$, so all three triangle inequalities follow from $\rho_f(p,q)\le2D$ and $D\le\rho_f(p,q)+D$.

Finally, if $p\in B_n$, $q\in B_m$, and $r\in B_\ell$, where $n,m,\ell$ are pairwise distinct, then \eqref{eq:external-distance} implies
\begin{align*}
 \rho_f(p,r)
 &=10\max\{a_n,a_\ell\}\\
 &\le10\max\{a_n,a_m,a_\ell\}\\
 &=\max\bigl\{\rho_f(p,q),\rho_f(q,r)\bigr\}\\
 &\le\rho_f(p,q)+\rho_f(q,r).
\end{align*}
Again, the same estimate applies to every permutation of $p,q,r$. Thus $\rho_f$ satisfies the triangle inequality.

Every nonzero distance is at least $3a_1$, so $M_f$ is uniformly discrete and hence complete. For a fixed $u\in V$, we have $\rho_f((1,u),(n,u))=10a_n\to\infty$ as $n\to\infty$, so $M_f$ is unbounded. Its cardinality is $|\Omega_I|=\aleph_0\cdot\kappa=\kappa$.
\end{proof}

\begin{lemma}\label{lem:band}
Let $f\in\R^I$ and $n\ge2$. If $p,q\in\Omega_I$ satisfy $\rho_f(p,q)\in(2a_n,6a_n)$, then $p,q\in B_n$. More precisely, if $p,q$ do not both belong to $B_n$, then
\begin{equation}\label{eq:band-gap}
 \rho_f(p,q)\le\frac{a_n}{10}
 \qquad\text{or}\qquad
 \rho_f(p,q)\ge10a_n.
\end{equation}
\end{lemma}

\begin{proof}
If $p=q$, then $\rho_f(p,q)=0$, so assume $p\ne q$. If both points belong to blocks of index less than $n$, then $\rho_f(p,q)\le10a_{n-1}=a_n/10$, whether or not they belong to the same block. If the points belong to different blocks and at least one index is at least $n$, then their distance is at least $10a_n$. The remaining possibility is that both points belong to the same block $B_m$ with $m>n$, in which case $\rho_f(p,q)\ge3a_m\ge3a_{n+1}=300a_n$.
\end{proof}

\begin{lemma}\label{lem:tail}
Let $f,g\in\R^I$, and suppose $R\in\Corr(M_f,M_g)$ has distortion $C:=\dis R<\infty$. Choose an integer $N\ge2$ such that $a_n>C+4n$ for all $n\ge N$. Then, for every $n\ge N$ and $u\in V$,
\begin{equation}\label{eq:tail-identity}
 R((n,u))=\{(n,u)\}.
\end{equation}
\end{lemma}

\begin{proof}
Such an $N$ exists because $a_n-4n\to\infty$. Fix $n\ge N$.

We first show that the $n$th block is preserved. For distinct $p,q\in B_n$, take $p'\in R(p)$ and $q'\in R(q)$. By \cref{prop:metric} and $C<a_n$,
\begin{equation}\label{eq:image-band}
 2a_n<3a_n-C\le \rho_g(p',q')\le5a_n+C<6a_n.
\end{equation}
By \cref{lem:band}, $p',q'\in B_n$. Since for every $p\in B_n$ one can choose a different point $q\in B_n$, this gives $R[B_n]\subseteq B_n$. Applying the same argument to the inverse correspondence $R^{-1}$ gives $R^{-1}[B_n]\subseteq B_n$. Consequently, $R$ contains no pair connecting $B_n$ to another block, and its restriction to $B_n\times B_n$ is a correspondence.

Next we show that this restricted correspondence is single-valued and injective. If one point $p\in B_n$ corresponded to two distinct points $p',p''\in B_n$, then
\[
 C\ge |\rho_f(p,p)-\rho_g(p',p'')|=\rho_g(p',p'')\ge3a_n>C,
\]
a contradiction. Similarly, applying the argument to the inverse correspondence shows that two different points cannot correspond to the same point. Thus the restricted correspondence is the graph of a bijection $\sigma_n:V\to V$.

Finally, suppose $\sigma_n\ne\id_V$. By \eqref{eq:adjacency-change}, there are distinct vertices $u,v$ for which the adjacency indicator changes. The corresponding difference of distances satisfies
\begin{align*}
 &\bigl|\rho_f((n,u),(n,v))
       -\rho_g((n,\sigma_nu),(n,\sigma_nv))\bigr|\\
 &\quad=\left|a_n\bigl(\varepsilon_H(u,v)-\varepsilon_H(\sigma_nu,\sigma_nv)\bigr)
       +b_{n,f}(u,v)-b_{n,g}(\sigma_nu,\sigma_nv)\right|\\
 &\quad\ge a_n-4n>C,
\end{align*}
contradicting the distortion bound. Therefore $\sigma_n=\id_V$, which proves \eqref{eq:tail-identity}.
\end{proof}

\begin{proof}[Proof of \cref{thm:coding}]
Fix $f,g\in\R^I$ and put $L:=\Delta_I(f,g)\in[0,\infty]$.

For the upper bound, use the identity correspondence $I_{\Omega_I}$. The distances between blocks do not depend on $f$, and the only distances that vary within a block are those of vertex pairs $\{o,i\}$ with $i\in I$. Thus
\begin{equation}\label{eq:id-distortion}
 \dis I_{\Omega_I}
 =\sup_{\omega,\omega'\in\Omega_I}|\rho_f(\omega,\omega')-\rho_g(\omega,\omega')|
 =2\sup_{\substack{n\ge1\\i\in I}}
 |c_n(f(i))-c_n(g(i))|.
\end{equation}
By \eqref{eq:clipping-properties},
\begin{equation}\label{eq:clipping-sup-equality}
 \sup_{\substack{n\ge1\\i\in I}}
 |c_n(f(i))-c_n(g(i))|
 =\sup_{i\in I}|f(i)-g(i)|=L.
\end{equation}
Indeed, $|c_n(f(i))-c_n(g(i))|\le|f(i)-g(i)|$, and for each $i\in I$ we may take an integer $n\ge\max\{1,|f(i)|,|g(i)|\}$, so that $c_n(f(i))=f(i)$ and $c_n(g(i))=g(i)$. It follows from \cref{prop:correspondence} that $\dGH(M_f,M_g)\le L$.

For the lower bound, take an arbitrary $R\in\Corr(M_f,M_g)$. If $\dis R=\infty$, there is nothing to prove. Otherwise, put $C:=\dis R<\infty$ and choose $N$ as in \cref{lem:tail}. Fix any $i\in I$ and choose an integer $n\ge\max\{N,|f(i)|,|g(i)|\}$. By \cref{lem:tail}, both $((n,o),(n,o))$ and $((n,i),(n,i))$ belong to $R$. Moreover, \eqref{eq:clipping-properties} gives $c_n(f(i))=f(i)$ and $c_n(g(i))=g(i)$, so
\[
 C\ge |\rho_f((n,o),(n,i))-\rho_g((n,o),(n,i))|
   =2|f(i)-g(i)|.
\]
As $i$ is arbitrary, $C\ge2L$. In particular, if $L=\infty$, no correspondence of finite distortion exists. If $L<\infty$, taking the infimum over all correspondences gives $\dGH(M_f,M_g)\ge L$. Together with the upper bound this proves \eqref{eq:coding}. When $L<\infty$, \eqref{eq:id-distortion} also shows that the distortion of the identity correspondence is $2L$, proving its optimality.
\end{proof}

\begin{definition}\label{def:affine-cloud}
For $f\in\R^I$, set $\kappa_I:=\max\{|I|,\aleph_0\}$ and define
\begin{equation}\label{eq:affine-cloud}
 E_f:\ell^\infty(I)\longrightarrow\GH_{\le\kappa_I},
 \qquad E_f(h):=[M_{f+h}].
\end{equation}
\end{definition}

By \cref{thm:coding}, for all $h,k\in\ell^\infty(I)$,
\begin{equation}\label{eq:affine-cloud-isometry}
 \dGH(M_{f+h},M_{f+k})=\|h-k\|_\infty.
\end{equation}
Taking $k=0$, we also have $\dGH(M_{f+h},M_f)=\|h\|_\infty<\infty$. Therefore $E_f$ gives an isometric embedding of $\ell^\infty(I)$ into $\mathcal C_{\le\kappa_I}(M_f)$.

\subsection{Proof of the main theorem}\label{sec:extended-universality}

\begin{lemma}\label{lem:extended-Kuratowski}
Let $(X,d)$ be a nonempty extended metric space with a nonempty dense subset $D\subseteq X$, and put $I:=\N_{\ge1}\times D$. For each $x\in X$, define
\begin{equation}\label{eq:truncated-coordinates}
 \varphi_x(m,z):=\min\{d(x,z),m\}
 \qquad((m,z)\in I),
\end{equation}
where $\min\{\infty,m\}=m$. Then $\varphi_x\in\R^I$, and
\begin{equation}\label{eq:extended-Kuratowski}
 \Delta_I(\varphi_x,\varphi_y)=d(x,y)
 \qquad(x,y\in X).
\end{equation}
\end{lemma}

\begin{proof}
First suppose $r:=d(x,y)<\infty$. If $z$ belongs to the same finite-distance component as $x$ and $y$, the reverse triangle inequality and the $1$-Lipschitz property of $t\mapsto\min\{t,m\}$ give $|\varphi_x(m,z)-\varphi_y(m,z)|\le r$. If $z$ belongs to a different component, both values equal $m$. Thus, in either case, $\Delta_I(\varphi_x,\varphi_y)\le r$.

Conversely, by \cref{lem:component-density}, $D$ meets each component densely. For every $\varepsilon>0$, we may therefore choose $z\in D$ with $d(x,z)<\varepsilon$. The point $z$ belongs to the same component as $x,y$, and
\[
 r-d(x,z)\le d(y,z)\le r+d(x,z)<r+\varepsilon.
\]
Choose an integer $m>r+\varepsilon$. Then $d(x,z)<m$ and $d(y,z)<m$, so $\varphi_x(m,z)=d(x,z)$ and $\varphi_y(m,z)=d(y,z)$. Hence
\[
 |\varphi_x(m,z)-\varphi_y(m,z)|
 =|d(x,z)-d(y,z)|
 \ge r-2d(x,z)>r-2\varepsilon.
\]
Letting $\varepsilon\downarrow0$ yields equality.

Next suppose $d(x,y)=\infty$. Choose $z\in D$ with $d(x,z)<1$. The triangle inequality implies $d(y,z)=\infty$. Thus, writing $a:=d(x,z)<1$, we obtain $|\varphi_x(m,z)-\varphi_y(m,z)|=m-a\to\infty$ as $m\to\infty$. This proves \eqref{eq:extended-Kuratowski} also for infinite distances.
\end{proof}

\begin{proof}[Proof of \cref{thm:main}]
Choose a nonempty dense subset $D\subseteq X$ with $|D|=\dens X$ and put $I:=\N_{\ge1}\times D$ and $\kappa:=\max\{\dens X,\aleph_0\}$. Apply \cref{thm:coding} to the functions $\varphi_x$ of \cref{lem:extended-Kuratowski}, and set $\Omega:=\Omega_I$ and $\rho_x:=\rho_{\varphi_x}$ for $x\in X$. Then the map
\begin{equation}\label{eq:embedding-X}
 \Phi:X\longrightarrow\GH_{\le\kappa},\qquad
 \Phi(x):=[(\Omega,\rho_x)]=[M_{\varphi_x}]
\end{equation}
satisfies
\[
 \dGH\bigl(\Phi(x),\Phi(y)\bigr)
 =\Delta_I(\varphi_x,\varphi_y)=d_X(x,y)
 \qquad(x,y\in X).
\]
Since the distance between distinct points is positive or infinite, $\Phi$ is also injective and hence is the required isometric embedding. The cardinality assertion in \eqref{eq:main} follows from $|I|=\max\{\dens X,\aleph_0\}$. The properties of the representing spaces and the optimality of the identity correspondence for pairs at finite distance follow from \cref{thm:coding}. Moreover, preservation of distances implies that $x,y$ belong to the same finite-distance component if and only if $\Phi(x),\Phi(y)$ belong to the same cloud.
\end{proof}

\subsection{Proper representatives in the countable case}\label{sec:proper}

The spaces $M_f$ are not proper, since each block is infinite, bounded, and uniformly discrete. For countable $I$, we replace these blocks by finite sets.

Let $I=\{i_1,i_2,\ldots\}$ be countably infinite, with the displayed enumeration fixed.

\begin{lemma}\label{lem:finite-rigid-trees}
For every integer $n\ge1$, there exist a finite simple undirected graph $H_n=(V_n,E_n)$ and distinct vertices $o_n,i_{1,n},\ldots,i_{n,n}\in V_n$ such that $H_n$ is a tree and $\Aut(H_n)=\{\id_{V_n}\}$.
\end{lemma}

\begin{proof}
Fix $n\ge1$. Take a vertex $v_n^*$ and set
\[
 V_n:=\{v_n^*\}\sqcup
 \bigl\{(k,r):1\le k\le n+3,\ 1\le r\le k\bigr\}.
\]
Define the edge set by
\[
 \begin{aligned}
 E_n={}&\bigl\{\{v_n^*,(k,1)\}:1\le k\le n+3\bigr\}\\
 &\cup\bigl\{\{(k,r),(k,r+1)\}:
      1\le k\le n+3,\ 1\le r<k\bigr\}.
 \end{aligned}
\]
The graph $H_n=(V_n,E_n)$ is connected and has no cycles, so it is a tree. The vertex $v_n^*$ has degree $n+3\ge4$, whereas every other vertex has degree at most $2$. Thus every automorphism fixes $v_n^*$. The components of $H_n-\{v_n^*\}$ are paths with distinct cardinalities $1,\ldots,n+3$, so each automorphism preserves every component. It also fixes the endpoint adjacent to $v_n^*$, and hence every vertex of that path. Therefore $\Aut(H_n)=\{\id_{V_n}\}$.

Finally, the choices $o_n:=(1,1)$ and $i_{j,n}:=(j+1,j+1)$, for $1\le j\le n$, give distinct vertices of degree $1$.
\end{proof}

Fix a tree and distinguished vertices as in \cref{lem:finite-rigid-trees} for each $n$.

\begin{definition}\label{def:proper-metric}
For distinct $u,v\in V_n$, let $\varepsilon_n(u,v)$ denote the adjacency indicator of $H_n$, defined as in \eqref{eq:adjacency}. Use $a_n=100^n$ and $c_n$ from \eqref{eq:clipping-scale}, and for $f\in\R^I$ put
\[
 \widehat b_{n,f}(u,v)
 :=\begin{cases}
  2\bigl(n+c_n(f(i_j))\bigr),
  &\text{if }\{u,v\}=\{o_n,i_{j,n}\}\text{ for some }1\le j\le n,\\
  0,&\text{otherwise}.
 \end{cases}
\]
Set $\widehat B_n:=\{n\}\times V_n$ and $\widehat\Omega_I:=\bigsqcup_{n\ge1}\widehat B_n$. For each $f\in\R^I$, define a symmetric function $\widehat\rho_f$ to be zero on the diagonal and, for distinct points, by
\begin{align}
 \widehat\rho_f((n,u),(n,v))
 &=3a_n+a_n\varepsilon_n(u,v)+\widehat b_{n,f}(u,v),
 &&u\ne v,
 \label{eq:proper-internal}\\
 \widehat\rho_f((n,u),(m,v))
 &=10\max\{a_n,a_m\},
 &&n\ne m.
 \label{eq:proper-external}
\end{align}
Write $\widehat M_f:=(\widehat\Omega_I,\widehat\rho_f)$.
\end{definition}

\begin{proposition}\label{prop:proper-metric}
For every $f\in\R^I$, $\widehat\rho_f$ is a finite-valued metric, and $3a_n\le\widehat\rho_f(p,q)\le5a_n$ for distinct $p,q\in\widehat B_n$. Moreover, $\widehat M_f$ is countable, complete, proper, uniformly discrete, and unbounded, and for every $R>0$,
\begin{equation}\label{eq:uniform-bounded-geometry}
 \sup_{f\in\R^I}\sup_{z\in\widehat\Omega_I}
 |\overline B_{\widehat\rho_f}(z,R)|<\infty.
\end{equation}
\end{proposition}

\begin{proof}
The bound $0\le\widehat b_{n,f}(u,v)\le4n\le a_n$ gives $3a_n\le\widehat\rho_f(p,q)\le5a_n$ for distinct $p,q\in\widehat B_n$. These bounds and \eqref{eq:proper-external} imply the triangle inequality by the same three cases as in the proof of \cref{prop:metric}; the remaining metric axioms follow directly from the definition.

Since each $V_n$ is finite and nonempty, $\widehat\Omega_I$ is countably infinite. Every nonzero distance is at least $3a_1$, so $\widehat M_f$ is uniformly discrete and complete. For $n\ge2$, the distance $\widehat\rho_f((1,o_1),(n,o_n))=10a_n$ tends to infinity with $n$, so $\widehat M_f$ is unbounded.

Fix $R>0$, choose an integer $N=N(R)\ge2$ such that $3a_n>R$ for $n\ge N$, and put $F_R:=\bigsqcup_{n<N}\widehat B_n$. The within-block lower bound and \eqref{eq:proper-external} give, for every $f\in\R^I$,
\[
 \begin{aligned}
 \overline B_{\widehat\rho_f}(z,R)&=\{z\},&&z\notin F_R,\\
 \overline B_{\widehat\rho_f}(z,R)&\subseteq F_R,&&z\in F_R.
 \end{aligned}
\]
Consequently,
\[
 \sup_{f\in\R^I}\sup_{z\in\widehat\Omega_I}
 |\overline B_{\widehat\rho_f}(z,R)|
 \le\max\left\{1,\sum_{n<N}|V_n|\right\}<\infty.
\]
This proves \eqref{eq:uniform-bounded-geometry}. Every closed bounded subset is contained in a finite closed ball and is therefore finite and compact. Hence $\widehat M_f$ is proper.
\end{proof}

\begin{lemma}\label{lem:proper-tail}
Let $f,g\in\R^I$ and suppose a correspondence $R\in\Corr(\widehat M_f,\widehat M_g)$ has distortion $C:=\dis R<\infty$. If $N\ge2$ is an integer such that $a_n>C+4n$ for all $n\ge N$, then $R((n,u))=\{(n,u)\}$ for all $n\ge N$ and $u\in V_n$.
\end{lemma}

\begin{proof}
Fix $n\ge N$. By \cref{prop:proper-metric} and \eqref{eq:proper-external}, the argument of \cref{lem:band} gives
\[
 \widehat\rho_h(p,q)\in(2a_n,6a_n)
 \quad\Longrightarrow\quad p,q\in\widehat B_n
 \qquad(h\in\{f,g\}).
\]
Since $C<a_n$ and $|V_n|\ge2$, the block-preservation argument in the proof of \cref{lem:tail}, applied to $R$ and $R^{-1}$, yields $R[\widehat B_n]\subseteq\widehat B_n$ and $R^{-1}[\widehat B_n]\subseteq\widehat B_n$. The lower bound $3a_n>C$ then excludes multiple images and multiple preimages, as in the same proof. Thus there is a bijection $\sigma_n:V_n\to V_n$ such that $R((n,u))=\{(n,\sigma_nu)\}$ for every $u\in V_n$.

If $\sigma_n\ne\id_{V_n}$, the rigidity of $H_n$ gives distinct $u,v\in V_n$ with $\varepsilon_n(u,v)\ne\varepsilon_n(\sigma_nu,\sigma_nv)$. Since both perturbation terms belong to $[0,4n]$, \eqref{eq:proper-internal} gives
\[
 C\ge
 \bigl|\widehat\rho_f((n,u),(n,v))
       -\widehat\rho_g((n,\sigma_nu),(n,\sigma_nv))\bigr|
 \ge a_n-4n>C,
\]
a contradiction. Hence $\sigma_n=\id_{V_n}$, which proves the assertion.
\end{proof}

\begin{theorem}\label{thm:proper-coding}
Let $I$ be a countably infinite set. For the family of metrics in \cref{def:proper-metric}, the map
\[
 \widehat{\mathcal E}_I:(\R^I,\Delta_I)\longrightarrow(\GH_{\le\aleph_0},\dGH),\qquad
 \widehat{\mathcal E}_I(f):=[\widehat M_f]
\]
is an isometric embedding. That is, for all $f,g\in\R^I$,
\begin{equation}\label{eq:proper-coding}
 \begin{aligned}
 \dGH(\widehat M_f,\widehat M_g)&=\Delta_I(f,g)\\
 &=\frac12\sup_{\omega,\omega'\in\widehat\Omega_I}
 |\widehat\rho_f(\omega,\omega')-\widehat\rho_g(\omega,\omega')|.
 \end{aligned}
\end{equation}
Each representative is a complete, proper, uniformly discrete, and unbounded metric space on the common countable underlying set $\widehat\Omega_I$, and satisfies \eqref{eq:uniform-bounded-geometry}. If $\Delta_I(f,g)<\infty$, the identity relation $I_{\widehat\Omega_I}$ is an optimal correspondence.
\end{theorem}

\begin{proof}
Fix $f,g\in\R^I$ and put $L:=\Delta_I(f,g)\in[0,\infty]$. The identity correspondence satisfies
\[
 \dis I_{\widehat\Omega_I}
 =2\sup_{n\ge1}\sup_{1\le j\le n}
   |c_n(f(i_j))-c_n(g(i_j))|
 =2L.
\]
For the last equality, \eqref{eq:clipping-properties} gives the upper bound, while for each $j\ge1$ an integer $n\ge\max\{j,|f(i_j)|,|g(i_j)|\}$ satisfies $c_n(f(i_j))=f(i_j)$ and $c_n(g(i_j))=g(i_j)$. By \cref{prop:correspondence}, $\dGH(\widehat M_f,\widehat M_g)\le L$.

For the reverse inequality, let $R\in\Corr(\widehat M_f,\widehat M_g)$ have finite distortion $C$, and choose $N$ as in \cref{lem:proper-tail}. Fix $j\ge1$ and take an integer $n\ge\max\{N,j,|f(i_j)|,|g(i_j)|\}$. By \cref{lem:proper-tail}, $R$ contains the pairs $((n,o_n),(n,o_n))$ and $((n,i_{j,n}),(n,i_{j,n}))$. Hence
\[
 C\ge
 \bigl|\widehat\rho_f((n,o_n),(n,i_{j,n}))
       -\widehat\rho_g((n,o_n),(n,i_{j,n}))\bigr|
 =2|f(i_j)-g(i_j)|.
\]
Taking the supremum over $j$ gives $C\ge2L$. If $L=\infty$, no finite-distortion correspondence exists, so \cref{prop:correspondence} gives $\dGH(\widehat M_f,\widehat M_g)=\infty$. If $L<\infty$, the same proposition gives the lower bound $\dGH(\widehat M_f,\widehat M_g)\ge L$. Together with the calculation of $\dis I_{\widehat\Omega_I}$, this proves \eqref{eq:proper-coding} and the optimality of the identity correspondence when $L<\infty$. Distance preservation implies injectivity, and the properties of the representatives follow from \cref{prop:proper-metric}.
\end{proof}

For a nonempty finite index set $I$, we can extend each function by zero along an injection $j:I\to\N_{\ge1}$ and apply \cref{thm:proper-coding} without changing the uniform extended metric. Hence the conclusion of that theorem holds for every nonempty countable set, including finite sets.

\section{Clouds and spheres}
\label{sec:clouds}

\subsection{Finite-distance components and universal clouds}

We use the following form of the Kuratowski embedding.

\begin{lemma}\label{lem:Kuratowski}
Let $(X,d_X)$ be a nonempty metric space, $x_0\in X$ a base point, and $D\subseteq X$ a nonempty dense subset. For each $x\in X$, define
\begin{equation}\label{eq:Kuratowski}
 K_x(z):=d_X(x,z)-d_X(x_0,z)
 \qquad(x\in X,\ z\in D).
\end{equation}
Then $K_x\in\ell^\infty(D)$ for every $x\in X$, and
\begin{equation}\label{eq:Kuratowski-isometry}
 \|K_x\|_\infty=d_X(x,x_0),
 \qquad
 \|K_x-K_y\|_\infty=d_X(x,y).
\end{equation}
Consequently, $x\mapsto K_x$ is an isometric embedding from $X$ into $\ell^\infty(D)$.
\end{lemma}

\begin{proof}
For every $x\in X$ and $z\in D$, the reverse triangle inequality gives $|K_x(z)|\le d_X(x,x_0)$. Hence $K_x\in\ell^\infty(D)$.

Fix $x,y\in X$. For every $z\in D$, we have $|K_x(z)-K_y(z)|=|d_X(x,z)-d_X(y,z)|\le d_X(x,y)$, so $\|K_x-K_y\|_\infty\le d_X(x,y)$. Conversely, let $\varepsilon>0$ and choose $z\in D$ with $d_X(x,z)<\varepsilon$. Then
\[
 \begin{aligned}
 \|K_x-K_y\|_\infty
 &\ge |d_X(x,z)-d_X(y,z)|\\
 &\ge d_X(x,y)-2d_X(x,z)
 >d_X(x,y)-2\varepsilon.
 \end{aligned}
\]
Letting $\varepsilon\downarrow0$ proves $\|K_x-K_y\|_\infty=d_X(x,y)$. Since $K_{x_0}=0$, setting $y=x_0$ also gives $\|K_x\|_\infty=d_X(x,x_0)$. The distance identity implies injectivity, so $x\mapsto K_x$ is an isometric embedding.
\end{proof}

\begin{theorem}\label{thm:universal-cloud}
For every infinite cardinal $\kappa$, there exists a complete, uniformly discrete, and unbounded metric space $U_\kappa$ of cardinality $\kappa$ such that every nonempty metric space $X$ with $\dens X\le\kappa$ admits an isometric embedding into $\mathcal C_{\le\kappa}(U_\kappa)$. Every point in the image can be represented by a complete, uniformly discrete, and unbounded metric space of cardinality $\kappa$.
\end{theorem}

\begin{proof}
Fix a set $I$ with $|I|=\kappa$ and the construction of \cref{def:multiscale-metric}, and set $U_\kappa:=M_0$, where $0$ is the zero function on $I$.

Let $(X,d_X)$ be a nonempty metric space with $\dens X\le\kappa$. Choose a nonempty dense subset $D\subseteq X$ with $|D|\le\kappa$, a point $x_0\in X$, and an injection $j:D\to I$. For the functions $K_x$ of \cref{lem:Kuratowski}, define
\begin{equation}\label{eq:padding}
 \widetilde K_x(i)
 :=\begin{cases}
 K_x(z),&i=j(z),\ z\in D,\\
 0,&i\notin j(D).
 \end{cases}
\end{equation}
Since $j$ is injective, $\widetilde K_x\in\ell^\infty(I)$ and extension by zero preserves the uniform norm. Thus \eqref{eq:affine-cloud-isometry} and \cref{lem:Kuratowski} give
\[
 \begin{aligned}
 \dGH(M_{\widetilde K_x},U_\kappa)
 &=\|\widetilde K_x\|_\infty
 =\|K_x\|_\infty=d_X(x,x_0)<\infty,\\
 \dGH(M_{\widetilde K_x},M_{\widetilde K_y})
 &=\|\widetilde K_x-\widetilde K_y\|_\infty
 =\|K_x-K_y\|_\infty=d_X(x,y)
 \qquad(x,y\in X).
 \end{aligned}
\]
By \cref{prop:metric}, $U_\kappa$ and every $M_{\widetilde K_x}$ are complete, uniformly discrete, and unbounded, with cardinality $\kappa$. Consequently, $x\mapsto[M_{\widetilde K_x}]$ is an isometric embedding into $\mathcal C_{\le\kappa}(U_\kappa)$ with the required representatives.
\end{proof}

\begin{definition}\label{def:universal-family}
Let $\kappa$ be an infinite cardinal, fix a set $J$ with $|J|=\kappa$, and put $I:=J\times\N_{\ge1}$. Fix the construction of \cref{def:multiscale-metric} for this set $I$, and set
\begin{equation}\label{eq:universal-family}
 \mathcal U_\kappa:=\{[M_f]\mid f\in\R^I\}
 \subseteq\GH_{\le\kappa},
\end{equation}
equipped with the restriction of $\dGH$. For $f\in\R^I$, write $\mathcal U_\kappa(f):=\{[M_{f+h}]\mid h\in\ell^\infty(I)\}$. Also, for each $A\subseteq J$, define $q_A\in\R^I$ by
\begin{equation}\label{eq:cloud-labels}
 q_A(j,m):=m\mathbf1_A(j)
 \qquad((j,m)\in I),
\end{equation}
where $\mathbf1_A(j)$ equals $1$ if $j\in A$ and $0$ otherwise.
\end{definition}

\begin{theorem}\label{thm:component-universality}
For every infinite cardinal $\kappa$, the family $\mathcal U_\kappa$ in \cref{def:universal-family} is an extended metric space with exactly $2^\kappa$ finite-distance components. These components are the sets $\mathcal U_\kappa(f)$, for $f\in\R^I$, and each is isometric to $\ell^\infty(\kappa)$.

A nonempty extended metric space $X$ admits an isometric embedding into $\mathcal U_\kappa$ if and only if it has at most $2^\kappa$ finite-distance components and each component admits an isometric embedding into $\ell^\infty(\kappa)$. In particular, every nonempty extended metric space with at most $2^\kappa$ finite-distance components, each of density character at most $\kappa$, admits an isometric embedding into $\mathcal U_\kappa$.

Every point of these embeddings can be represented by a complete, uniformly discrete, and unbounded metric space on a common underlying set of cardinality $\kappa$. Furthermore, for each $f\in\R^I$, the cardinality-restricted cloud $\mathcal C_{\le\kappa}(M_f)$ is isometrically universal for all metric spaces of density character at most $\kappa$, and there are $2^\kappa$ distinct cardinality-restricted clouds obtained in this way.
\end{theorem}

\begin{proof}
Use the sets $J,I$ and the family $(M_f)_{f\in\R^I}$ fixed in \cref{def:universal-family}. Define $\mathcal E:\R^I\to\mathcal U_\kappa$ by $\mathcal E(f):=[M_f]$. By \cref{thm:coding}, for all $f,g\in\R^I$,
\begin{equation}\label{eq:component-proof-isometry}
 \dGH(M_f,M_g)
 =\sup_{(j,m)\in J\times\N_{\ge1}}|f(j,m)-g(j,m)|
 =\Delta_I(f,g).
\end{equation}
The map $\mathcal E$ is surjective by definition and injective because $\Delta_I$ is an extended metric. Hence it is an isometry onto $\mathcal U_\kappa$. Moreover,
\[
 \dGH(M_f,M_g)<\infty
 \iff g-f\in\ell^\infty(I)
 \iff [M_g]\in\mathcal U_\kappa(f).
\]
Thus the finite-distance components are precisely the sets $\mathcal U_\kappa(f)$. The map $E_f:h\mapsto[M_{f+h}]$ is surjective onto $\mathcal U_\kappa(f)$ and is an isometry by \eqref{eq:affine-cloud-isometry}. Since $|I|=\kappa$, composition with a bijection $\kappa\to I$ identifies $\ell^\infty(I)$ isometrically with $\ell^\infty(\kappa)$.

For distinct $A,B\subseteq J$, choose $j\in A\triangle B$. Then $|q_A(j,m)-q_B(j,m)|=m$ for every $m\ge1$, so $\Delta_I(q_A,q_B)=\infty$. The cosets $q_A+\ell^\infty(I)$ are therefore pairwise distinct. It follows that
\[
 2^\kappa
 \le |\R^I/\ell^\infty(I)|
 \le |\R^I|
 =(2^{\aleph_0})^\kappa
 =2^\kappa.
\]
Each coset corresponds to one finite-distance component, so there are exactly $2^\kappa$ such components.

We next prove the necessary condition. Suppose $F:X\to\mathcal U_\kappa$ is an isometric embedding, and write $(X_\alpha)_{\alpha\in\mathcal A(X)}$ for the finite-distance components of $X$. Choose $x_\alpha\in X_\alpha$ and set $g_\alpha:=\mathcal E^{-1}(F(x_\alpha))$. For $x\in X_\alpha$, define
\[
 T_\alpha(x):=\mathcal E^{-1}(F(x))-g_\alpha.
\]
By \eqref{eq:component-proof-isometry}, $T_\alpha(x)\in\ell^\infty(I)$ because $d_X(x,x_\alpha)<\infty$, and
\[
 \|T_\alpha(x)-T_\alpha(y)\|_\infty
 =\dGH(F(x),F(y))=d_X(x,y)
 \qquad(x,y\in X_\alpha).
\]
Thus each $X_\alpha$ embeds isometrically into $\ell^\infty(\kappa)$. If $\alpha\ne\beta$, then $\Delta_I(g_\alpha,g_\beta)=d_X(x_\alpha,x_\beta)=\infty$. Hence $\alpha\mapsto g_\alpha+\ell^\infty(I)$ is injective, and $|\mathcal A(X)|\le2^\kappa$.

Conversely, suppose $|\mathcal A(X)|\le2^\kappa$ and each $X_\alpha$ admits an isometric embedding $H_\alpha:X_\alpha\to\ell^\infty(J)$. Choose pairwise distinct subsets $A_\alpha\subseteq J$, one for each $\alpha\in\mathcal A(X)$. For $x\in X_\alpha$, put $h_x:=H_\alpha(x)$ and define $f_x\in\R^I$ by
\begin{equation}\label{eq:component-code}
 f_x(j,m):=m\mathbf1_{A_\alpha}(j)+h_x(j)
 \qquad((j,m)\in I).
\end{equation}
If $x,y\in X_\alpha$, then
\[
 \Delta_I(f_x,f_y)
 =\sup_{j\in J}|h_x(j)-h_y(j)|=d_X(x,y).
\]
If $x\in X_\alpha$ and $y\in X_\beta$ with $\alpha\ne\beta$, choose $j\in A_\alpha\triangle A_\beta$. For every $m\ge1$,
\[
 |f_x(j,m)-f_y(j,m)|
 \ge m-|h_x(j)-h_y(j)|\longrightarrow\infty
 \qquad(m\to\infty).
\]
Thus $\Delta_I(f_x,f_y)=\infty=d_X(x,y)$ in this case. By \eqref{eq:component-proof-isometry}, $x\mapsto\mathcal E(f_x)$ is an isometric embedding of $X$ into $\mathcal U_\kappa$. All its representatives $M_{f_x}$ have the same underlying set $\Omega_I$, of cardinality $\kappa$, and the properties asserted in \cref{thm:coding}.

Suppose now that $\dens X_\alpha\le\kappa$ for every $\alpha$. Choose $x_\alpha\in X_\alpha$, a nonempty dense subset $D_\alpha\subseteq X_\alpha$ with $|D_\alpha|\le\kappa$, and an injection $j_\alpha:D_\alpha\to J$. For $x\in X_\alpha$, define
\[
 H_\alpha(x)(j):=
 \begin{cases}
 d_X(x,z)-d_X(x_\alpha,z),&j=j_\alpha(z),\ z\in D_\alpha,\\
 0,&j\notin j_\alpha(D_\alpha).
 \end{cases}
\]
By \cref{lem:Kuratowski}, $H_\alpha(x)\in\ell^\infty(J)$ and $\|H_\alpha(x)-H_\alpha(y)\|_\infty=d_X(x,y)$ for $x,y\in X_\alpha$. The sufficient condition just proved therefore applies.

For any nonempty metric space $(Z,d_Z)$ with $\dens Z\le\kappa$, the same construction with $I$ in place of $J$ gives an isometric embedding $H:Z\to\ell^\infty(I)$. Fix $f\in\R^I$. By \eqref{eq:affine-cloud-isometry},
\[
 \begin{aligned}
 \dGH(M_{f+H(z)},M_f)&=\|H(z)\|_\infty<\infty,\\
 \dGH(M_{f+H(z)},M_{f+H(w)})&=\|H(z)-H(w)\|_\infty=d_Z(z,w)
 \qquad(z,w\in Z).
 \end{aligned}
\]
Since these representatives have cardinality $\kappa$, $z\mapsto[M_{f+H(z)}]$ is an isometric embedding into $\mathcal C_{\le\kappa}(M_f)$.

Finally, for $f,g\in\R^I$,
\[
 \mathcal C_{\le\kappa}(M_f)=\mathcal C_{\le\kappa}(M_g)
 \iff \dGH(M_f,M_g)<\infty
 \iff g-f\in\ell^\infty(I).
\]
Indeed, finite distance implies equality of the two clouds by the triangle inequality; conversely, equality implies $[M_g]\in\mathcal C_{\le\kappa}(M_f)$ because $[M_g]\in\mathcal C_{\le\kappa}(M_g)$. Thus these cardinality-restricted clouds are in bijection with $\R^I/\ell^\infty(I)$, and their number is $2^\kappa$.
\end{proof}

\begin{proposition}\label{prop:cloud-count}
For an infinite cardinal $\kappa$, define an equivalence relation on $\GH_{\le\kappa}$ by $[A]\asymp[B]\iff\dGH(A,B)<\infty$. Then
\begin{equation}\label{eq:cloud-count}
 |\GH_{\le\kappa}|
 =\left|\GH_{\le\kappa}/\!\asymp\right|=2^\kappa.
\end{equation}
Of these equivalence classes, $2^\kappa$ are each isometrically universal for all metric spaces of density character at most $\kappa$.
\end{proposition}

\begin{proof}
Fix $J$, $I$, and $(M_f)_{f\in\R^I}$ as in \cref{def:universal-family}. For each $A\subseteq J$, the $\asymp$-equivalence class of $[M_{q_A}]$ is $\mathcal C_{\le\kappa}(M_{q_A})$. If $A\ne B$, choose $j\in A\triangle B$. By \eqref{eq:cloud-labels}, $|q_A(j,m)-q_B(j,m)|=m$ for every $m\ge1$, so \cref{thm:coding} gives
\[
 \dGH(M_{q_A},M_{q_B})=\Delta_I(q_A,q_B)=\infty.
\]
These equivalence classes are therefore pairwise distinct. Together with \cref{prop:GH-set-representation}, this gives
\[
 2^\kappa\le\left|\GH_{\le\kappa}/\!\asymp\right|
 \le|\GH_{\le\kappa}|\le2^\kappa.
\]
Hence both cardinalities equal $2^\kappa$. By \cref{thm:component-universality}, each of the $2^\kappa$ distinct classes $\mathcal C_{\le\kappa}(M_{q_A})$ is isometrically universal for metric spaces of density character at most $\kappa$.
\end{proof}

\begin{corollary}\label{cor:proper-componentwise}
Let $(X,d_X)$ be a nonempty extended metric space all of whose finite-distance components are separable. Then the following three conditions are equivalent.
\begin{enumerate}
\item The number of finite-distance components of $X$ is at most $2^{\aleph_0}$.
\item $X$ admits an isometric embedding into $\GH_{\le\aleph_0}$.
\item There exist a common countable set $\Omega$ and a family of metrics $(\rho_x)_{x\in X}$ on $\Omega$ such that every $(\Omega,\rho_x)$ is proper, uniformly discrete, unbounded, and of bounded geometry, and
\[
 \dGH\bigl((\Omega,\rho_x),(\Omega,\rho_y)\bigr)=d_X(x,y)
 \qquad(x,y\in X).
\]
The bounded-geometry estimates can be chosen independently of $x\in X$, and if $d_X(x,y)<\infty$, then $I_\Omega$ is an optimal correspondence between $(\Omega,\rho_x)$ and $(\Omega,\rho_y)$.
\end{enumerate}
\end{corollary}

\begin{proof}
We first prove $(1)\Rightarrow(3)$. Take $J=\N_{\ge1}$ and $I=J\times\N_{\ge1}$ in \cref{def:universal-family}. Since every finite-distance component of $X$ is separable, condition~$(1)$ and \cref{thm:component-universality} give an isometric embedding $F:X\to\mathcal U_{\aleph_0}$. By \cref{thm:coding} and the definition of $\mathcal U_{\aleph_0}$, for each $x\in X$ there is a unique $f_x\in\R^I$ such that $F(x)=[M_{f_x}]$. Hence
\[
 \Delta_I(f_x,f_y)=\dGH(F(x),F(y))=d_X(x,y)
 \qquad(x,y\in X).
\]
The set $I$ is countably infinite, so \cref{thm:proper-coding} applies. Put $\Omega:=\widehat\Omega_I$ and $\rho_x:=\widehat\rho_{f_x}$. Then $\Omega$ is countable and
\[
 \dGH\bigl((\Omega,\rho_x),(\Omega,\rho_y)\bigr)
 =\Delta_I(f_x,f_y)=d_X(x,y).
\]
Each $(\Omega,\rho_x)$ is proper, uniformly discrete, and unbounded. For every $R>0$, \cref{prop:proper-metric} gives
\[
 \begin{aligned}
 \sup_{x\in X}\sup_{\omega\in\Omega}
 |\overline B_{\rho_x}(\omega,R)|
 &\le\sup_{f\in\R^I}\sup_{\omega\in\widehat\Omega_I}
 |\overline B_{\widehat\rho_f}(\omega,R)|
 <\infty.
 \end{aligned}
\]
Thus the bounded-geometry estimate is independent of $x$. If $d_X(x,y)<\infty$, \cref{thm:proper-coding} also gives $\dis I_\Omega=2\Delta_I(f_x,f_y)=2d_X(x,y)$, so $I_\Omega$ is optimal.

For $(3)\Rightarrow(2)$, define $F(x):=[(\Omega,\rho_x)]$. Since $\Omega$ is countable, $F(x)\in\GH_{\le\aleph_0}$, and condition~$(3)$ gives $\dGH(F(x),F(y))=d_X(x,y)$. This equality also implies injectivity, so $F$ is an isometric embedding.

Finally, assume $(2)$ and let $F:X\to\GH_{\le\aleph_0}$ be an isometric embedding. Choose $x_\alpha\in X_\alpha$ for each $\alpha\in\mathcal A(X)$. If $\alpha\ne\beta$, then $\dGH(F(x_\alpha),F(x_\beta))=d_X(x_\alpha,x_\beta)=\infty$, so these images lie in distinct $\asymp$-equivalence classes. By \cref{prop:cloud-count},
\[
 |\mathcal A(X)|
 \le\left|\GH_{\le\aleph_0}/\!\asymp\right|
 =2^{\aleph_0},
\]
which is condition~$(1)$.
\end{proof}

\subsection{Isometric universality on spheres}\label{sec:spheres}

\begin{definition}\label{def:universal-subspace}
Let $\kappa$ be an infinite cardinal and $I$ a set with $|I|=\kappa$, and fix the construction of \cref{def:multiscale-metric}. For $f\in\R^I$, set
\[
 Y_f:=\mathcal Q_\kappa(M_f),\qquad
 \overline E_f:\ell^\infty(I)\longrightarrow Y_f,\qquad
 \overline E_f(h):=[M_{f+h}]_0,
\]
and write $\mathcal L_f:=\overline E_f(\ell^\infty(I))$ for its image. Here $[M]_0$ denotes the equivalence class of $[M]\in\mathcal C_{\le\kappa}(M_f)$ in the zero-distance quotient.

When $I$ is countably infinite, for the construction in \cref{def:proper-metric} also put
\[
 \widehat Y_f:=\mathcal Q_{\aleph_0}(\widehat M_f),\qquad
 \widehat E_f:\ell^\infty(I)\longrightarrow\widehat Y_f,\qquad
 \widehat E_f(h):=[\widehat M_{f+h}]_0,
\]
and define its image by $\widehat{\mathcal L}_f:=\widehat E_f(\ell^\infty(I))$.
\end{definition}

By \cref{thm:coding} and the definition of the metric on the zero-distance quotient, for all $h,k\in\ell^\infty(I)$,
\begin{equation}\label{eq:universal-subspace-isometry}
 \dGH(\overline E_f(h),\overline E_f(k))
 =\dGH(M_{f+h},M_{f+k})=\|h-k\|_\infty.
\end{equation}
Thus $\overline E_f$ is a surjective isometry from $\ell^\infty(I)$ onto $\mathcal L_f$. The same equality holds for $\widehat E_f$ in the countable case, by \cref{thm:proper-coding}.

\begin{definition}\label{def:doubling}
A metric space $(Z,d_Z)$ is \emph{doubling} if there exists an integer $N\ge1$ such that, for every $z\in Z$ and $R>0$, the open ball $B_{d_Z}(z,R)$ can be covered by at most $N$ balls of radius $R/2$. The space $Z$ is \emph{locally doubling} at $p\in Z$ if $p$ has a neighborhood that is doubling with respect to the induced metric.
\end{definition}

In the following theorem, balls and spheres for $\dGH$ are taken in $Y_f$, with the analogous convention for $\widehat Y_f$ in the countable case.

\begin{theorem}\label{thm:spherical-universality}
Let $\kappa,I,f$ be as in \cref{def:universal-subspace}. For every $p\in\mathcal L_f$ and $r>0$, every nonempty metric space $X$ satisfying $\dens X\le\kappa$ and $\diam X\le2r$ admits an isometric embedding
\[
 J:X\longrightarrow S_{\dGH}(p,r)\cap\mathcal L_f.
\]
Moreover, both $S_{\dGH}(p,r)\cap\mathcal L_f$ and $S_{\dGH}(p,r)$ have cardinality and density character equal to $2^\kappa$.

For every $R>0$ and every nonempty bounded metric space $X$ with $\dens X\le\kappa$, there exists $\lambda>0$ such that the metric space $\lambda X$, obtained by multiplying the metric by $\lambda$, admits an isometric embedding into $B_{\dGH}(p,R)\cap\mathcal L_f$. Moreover, $Y_f$ is neither locally compact nor locally doubling at any point of $\mathcal L_f$.

If $I$ is countably infinite, all these conclusions also hold with $Y_f$ and $\mathcal L_f$ replaced by $\widehat Y_f$ and $\widehat{\mathcal L}_f$, respectively. In that case, the points in the images of the embeddings can be represented by complete, proper, uniformly discrete, and unbounded metric spaces on a common countable underlying set, with bounded-geometry estimates independent of the image point.
\end{theorem}

\begin{proof}
Fix $p\in\mathcal L_f$ and write $p=\overline E_f(h_0)$ with $h_0\in\ell^\infty(I)$. Let $r>0$ and let $(X,d_X)$ be a nonempty metric space with $\dens X\le\kappa$ and $\diam X\le2r$. Choose a nonempty dense subset $D\subseteq X$ with $|D|\le\kappa$ and an injection $\iota:D\to I$. For $x\in X$, define
\[
 u_x(i):=
 \begin{cases}
 h_0(i)+d_X(x,z)-r,&i=\iota(z),\ z\in D,\\
 h_0(i),&i\notin\iota(D).
 \end{cases}
\]
Since $0\le d_X(x,z)\le2r$, we have $|u_x(i)-h_0(i)|\le r$ for every $i\in I$. Thus $u_x\in\ell^\infty(I)$ and $\|u_x-h_0\|_\infty\le r$. For any $\varepsilon>0$, choose $z\in D$ with $d_X(x,z)<\varepsilon$. Then
\[
 \|u_x-h_0\|_\infty
 \ge|d_X(x,z)-r|
 \ge r-d_X(x,z)>r-\varepsilon.
\]
Letting $\varepsilon\downarrow0$ gives $\|u_x-h_0\|_\infty=r$. Define $J(x):=\overline E_f(u_x)$. By \eqref{eq:universal-subspace-isometry} and \cref{lem:Kuratowski},
\[
 \begin{aligned}
 \dGH(J(x),p)&=\|u_x-h_0\|_\infty=r,\\
 \dGH(J(x),J(y))
 &=\|u_x-u_y\|_\infty
 =\sup_{z\in D}|d_X(x,z)-d_X(y,z)|
 =d_X(x,y)
 \qquad(x,y\in X).
 \end{aligned}
\]
Hence $J$ is an isometric embedding into $S_{\dGH}(p,r)\cap\mathcal L_f$.

For $t>0$ and $\eta\in\{-1,1\}^I$, put $P_{t,\eta}:=\overline E_f(h_0+t\eta)$ and $A_t:=\{P_{t,\eta}:\eta\in\{-1,1\}^I\}$. For distinct $\eta,\theta\in\{-1,1\}^I$, \eqref{eq:universal-subspace-isometry} gives
\[
 \dGH(P_{t,\eta},p)=t,\qquad
 \dGH(P_{t,\eta},P_{t,\theta})
 =t\|\eta-\theta\|_\infty=2t.
\]
Consequently, $A_t\subseteq S_{\dGH}(p,t)\cap\mathcal L_f$ and $|A_t|=2^\kappa$.

Let $Z$ be either $S_{\dGH}(p,r)\cap\mathcal L_f$ or $S_{\dGH}(p,r)$. For every $\eta\in\{-1,1\}^I$, the set $Z\cap B_{\dGH}(P_{r,\eta},r)$ is nonempty and open in $Z$. These sets are pairwise disjoint: a point $q$ in two of them would satisfy
\[
 2r=\dGH(P_{r,\eta},P_{r,\theta})
 \le\dGH(P_{r,\eta},q)+\dGH(q,P_{r,\theta})<2r
 \qquad(\eta\ne\theta).
\]
Every dense subset of $Z$ must meet each of these $2^\kappa$ disjoint sets. Using \cref{prop:GH-set-representation}, we obtain
\[
 2^\kappa\le\dens Z\le|Z|\le|Y_f|
 \le|\GH_{\le\kappa}|\le2^\kappa.
\]
This proves the cardinality and density assertions for both spheres.

Now let $R>0$ and let $X$ be a nonempty bounded metric space with $\dens X\le\kappa$. Set $\lambda:=R/(1+\diam X)>0$ and $r:=R/2$. Then $\dens(\lambda X)=\dens X$ and $\diam(\lambda X)=\lambda\diam X<R=2r$. The sphere result gives an isometric embedding of $\lambda X$ into $S_{\dGH}(p,r)\cap\mathcal L_f$, which is contained in $B_{\dGH}(p,R)\cap\mathcal L_f$ because $r<R$.

Let $W$ be any neighborhood of $p$ in $Y_f$, and put $d_W:=\dGH|_{W\times W}$. Choose $s>0$ such that $B_{\dGH}(p,s)\subseteq W$ and set $t:=3s/4$. Then
\[
 A_t\subseteq B_{d_W}(p,s)=B_{\dGH}(p,s),
 \qquad d_W(P,Q)=2t=3s/2
 \quad(P,Q\in A_t,\ P\ne Q).
\]
For $a\in W$, the ball $B_{d_W}(a,s/2)$ contains at most one point of $A_t$, since otherwise
\[
 3s/2=d_W(P,Q)
 \le d_W(P,a)+d_W(a,Q)<s.
\]
Thus $B_{d_W}(p,s)$ cannot be covered by finitely many balls of radius $s/2$, and $(W,d_W)$ is not doubling. The same argument shows that $W$ is not totally bounded and hence not compact. Since $W$ was arbitrary, $Y_f$ is neither locally doubling nor locally compact at $p$.

When $I$ is countably infinite, \cref{thm:proper-coding} and \cref{prop:GH-set-representation} give
\[
 \dGH(\widehat E_f(h),\widehat E_f(k))=\|h-k\|_\infty
 \quad(h,k\in\ell^\infty(I)),\qquad
 |\widehat Y_f|\le2^{\aleph_0}.
\]
The preceding constructions and estimates therefore apply with $\overline E_f$, $Y_f$, and $\mathcal L_f$ replaced by $\widehat E_f$, $\widehat Y_f$, and $\widehat{\mathcal L}_f$. Every point $\widehat E_f(h)$ is represented by $\widehat M_{f+h}$ on the common countable set $\widehat\Omega_I$. By \cref{prop:proper-metric}, these spaces are complete, proper, uniformly discrete, and unbounded, and \eqref{eq:uniform-bounded-geometry} gives the required bound independently of $h\in\ell^\infty(I)$.
\end{proof}

\section*{Acknowledgments}
The author used GPT-6.0 models as AI-assisted tools in preparing this manuscript. The author reviewed and revised the mathematical content and takes full responsibility for the final manuscript.

\bigskip
\noindent\textsc{Mathematical Institute, Tohoku University}\\
Sendai 980-8578, Japan\\
\textit{Email address:} \texttt{fukano.ryo.t5@dc.tohoku.ac.jp}

\end{document}